\documentclass[final]{elsart1p}

\usepackage[dvips]{epsfig} 
\usepackage[latin1]{inputenc}
\usepackage{amsmath}
\usepackage{amsfonts}
\usepackage{amssymb}
\usepackage{enumitem}

\newcommand{\bp}[1]{\big(#1\big)} 
\newcommand{\Bp}[1]{\Big(#1\Big)}

\usepackage{lineno}

\begin{document}
\begin{frontmatter}


\title{Extensions and applications of ACF mappings}
\author{Jean-Philippe Chancelier}
\ead{jpc@cermics.enpc.fr}

\address{Université Paris-Est, CERMICS (ENPC), 6-8 Avenue Blaise Pascal,
Cité Descartes , F-77455 Marne-la-Vallée}



\begin{abstract}
Using a definition of ASF sequences derived from the definition of 
asymptotic contractions of the final type of ACF, we give some new 
fixed points theorem for cyclic mappings and  alternating mapping which 
extend results from \cite[Theorem 2]{suzuki2007} and \cite[Theorem 1]{zhang2007}. 
\end{abstract} 

\end{frontmatter}

\def\Fix{\mathop{\normalfont Fix}}
\def\defpar{\stackrel{\mbox{\tiny def}}{=}}
\def\argmax{\mathop{\mbox{\rm Argmax}}}
\def\argmin{\mathop{\mbox{\rm Argmin}}}
\def\bbP{{\mathbb P}} 
\def\bbC{{\mathbb C}} 
\def\bbE{{\mathbb E}} 
\def\E{{\cal E}}
\def\F{{\cal F}}
\def\H{{\cal H}}
\def\V{{\cal V}}
\def\W{{\cal W}}
\def\U{{\cal U}}
\def\R{{\mathbb R}}
\def\RB{{\mathbb R}}
\def\N{{\mathbb N}}
\def\M{{\mathbb M}}
\def\bbR{{\mathbb R}} 
\def\bbI{{\mathbb I}} 
\def\bbU{{\mathbb U}} 
\def\Tad{{{\cal T}_{\mbox{\tiny ad}}}}
\def\texte#1{\quad\mbox{#1}\quad}
\def\Proba#1{\bbP\left\{ #1 \right\}} 
\def\Probax#1#2{{\bbP}_{#1}\left\{ #2 \right\}} 
\def\ProbaU#1#2{{\bbP}^{#1} \left\{ #2 \right\}} 
\def\ProbaxU#1#2#3{{\bbP}^{#1}_{#2} \left\{ #3 \right\}} 
\def\valmoy#1{\bbE\left[ #1 \right]}
\def\valmoyDebut#1{\bbE [ #1 } 
\def\valmoyFin#1{ #1 ]} 
\def\valmoyp#1#2{\bbE_{#1}\left[ #2 \right]}
\def\valmoypDebut#1#2{\bbE_{#1} \left[ #2 \right.} 
\def\valmoypFin#1{ \left. #1 \right]} 
\def\valmoypU#1#2#3{\bbE_{#1}^{#2}\left[ #3 \right]}
\def\norminf#1{ {\Vert #1 \Vert}_{\infty}}
\def\norm#1{ {\Vert #1 \Vert}}
\def\Hun{${\text{\bf H}}_1$}
\def\Hdeux{${\text{\bf H}}_2$}
\def\Htrois{${\text{\bf H}}_3$}
\def\psca#1{\left< #1 \right>}
\def\slim{\sigma\mbox{-}\lim}
\def\seq#1{\{{#1}_n\}_{n\in\N}} 
\def\X{{\mathbb X}}


\newenvironment{myproof}{{\vspace{0.5cm}\small{\it Proof:}}}{\hfill$\Box$\normalsize
\\\smallskip}

\setenumerate{labelindent=\parindent,label=\emph{($\mbox{C}_{\arabic*}$)},ref=($\mbox{C}_{\arabic*}$)}

\section{Introduction} 

Many extensions of the well known Banach contraction principle 
\cite{banach} have been proposed in nonlinear analysis literature.
Among them fixed point theorems for Meir-Keeler contraction have been 
extensively studied \cite{meir-keeler,kirk-2,suzuki-2006} and a final (in some 
sense) generalization defined as \emph{asymptotic contraction of the final type} 
(ACF, for short) has been stated by T.Suzuki \cite[Theorem 5]{suzuki2007}.  
Our aim in this paper is to extend the results of T.Suzuki to more general cases 
with regards to the mappings. More precisely, we want to be able to use 
the same framework for proving fixed point theorems for alternating mappings 
\S\ref{alternate} or for cyclic mappings \S\ref{cycling}. For that purpose 
we propose the definition of $p$-ASF-$1$ and $p$-ASF-$2$ sequences which are defined 
without references to a mapping and prove some Cauchy properties of such sequences 
in Theorem~\ref{thmpasf}. In \S\ref{defacf}, we recall the definition of ACF mapping 
and relate ACF mapping to $p$-ASF mappings. When the $p$-ASF sequences are generated 
using $\{T^nx\}$ we show that the two definitions coincide (Theorem~\ref{asf-acf}).
We give an application to cyclic mappings in \S\ref{cycling} by providing a fixed 
point theorem which extends \cite[Theorem 2]{suzuki2007} to continuous $p$-ASF mappings. 
In \S\ref{alternate} we give an application to alternating mapping through
Theorem~\ref{thm:ptfixe} which extends the results of~\cite{zhang2007}.

\section{ACF sequences}

In \cite{suzuki2007} T.Suzuki introduces the definition of an {\em asymptotic
  contraction of the final type} (ACF, for short) and proves that if a mapping
$T$ is ACF then the sequence $\seq{x}$ defined by $x_n\defpar T^n x$ is a Cauchy
sequence for all $x\in \X$. Since our aim is to extend T.Suzuki results when
sequences $\seq{x}$ are generated by more general processes, we introduce a new
definition that we call ASF, which stands for {\em asymptotic sequences of the
  final type}.  The definition characterizes two sequences and not a
mapping. The link between the two definitions is the following. Suppose that the
mapping $T$ is ACF and for $x$, $y \in \X$ define two sequences $\seq{x}$,
$\seq{y}$ by $x_n\defpar T^n x$ and $y_n\defpar T^n y$. If for all $n\in \N$ we
have $x_n \not=y_n$ then the two sequences are ASF.  Properties of ASF sequences
are given in Lemma~\ref{lem:ASF} and a proof is given but note that the proof is
mostly a simple rephrase of \cite[Lemma 1 and 2]{suzuki2007}.  We first start by
the ASF definition.

In the sequel $(\X,d)$ is a complete metric space and $p$ is a given function 
from $\X\times \X$ into $[0,\infty)$.

\begin{defn} \label{def:acf} We say that two sequences $\seq{x}$, $\seq{y}$ with
  $x_n,y_n \in \X$ are $p$-ASF-$1$ if the following are satisfied:
  \begin{enumerate}
  \item \label{it:acfun} For each $\epsilon >0$ there exists $\delta >0$ such that 
    if for $i \in \N$ we have $p(x_i,y_i) < \delta$ then 
    $\limsup_{n \to \infty} p(x_{n},y_{n}) \le \epsilon$\,;
  \item \label{it:acfdeux} For each $\epsilon > 0$, there exists $\delta > 0$ 
    such that for $i$,$j \in \N$ with $\epsilon < p(x_i, y_i) < \epsilon + \delta$, 
    there exists $\nu \in \N$ such that $p(x_{\nu+i},y_{\nu+i}) \le \epsilon$\,;
  \item \label{it:acftrois} For each given $(x_i,y_i)$ such that $p(x_i,y_i)\ne 0$ 
    there exists $\nu \in \N$ such that $$p(x_{\nu+i}, y_{\nu+i}) < p(x_i,y_i)\,.$$
  \end{enumerate}
\end{defn}

\begin{lem}\label{lem:ASF} Let  $\seq{x}$, $\seq{y}$ be two $p$-ASF-1 sequences 
then $\lim_{n\to \infty} p(x_n,y_n)=0$. 
\end{lem}

\begin{myproof}  
  We follow \cite[Lemma 2]{suzuki2007}.
  If we suppose that there exists $i\in\N$ such that $p(x_i,y_i)=0$, 
  we conclude directly using \ref{it:acfun} that $\lim_{n\to \infty} p(x_n,y_n)=0$. 
  Thus we assume now that $p(x_n,y_n)\not=0$ for all $n\in \N$.
  We first prove that if $\seq{x}$, $\seq{y}$ satisfy \ref{it:acfdeux} and
  \ref{it:acftrois} then $\liminf_{n \to \infty} p(x_n,y_n)=0$.
  Using the fact that $p$ is nonnegative and repeatedly using Property~\ref{it:acftrois} 
  it is possible to build an extracted decreasing sub-sequence 
  $p(x_{\sigma(n)},y_{\sigma(n)})$ such that 
  $0 \le p(x_{\sigma(n)},y_{\sigma(n)}) \le p(x_0,y_0)$ which implies that 
  $\liminf_{n \to \infty} p(x_n,y_n)=\alpha$ exists and is finite. 
  Suppose that $\alpha >0$. We first show that we must have $\alpha <
  p(x_n,y_n)$ for all $n \in \N$.  Indeed suppose that there exists $n_0$ such
  that $p(x_{n_0},y_{n_0}) \le \alpha$ then repeatedly using 
  \ref{it:acftrois} we
  can build an extracted decreasing sequence $p(x_{\sigma(n)},y_{\sigma(n)})$ such
  that $p(x_{\sigma(n)},y_{\sigma(n)}) < p(x_{n_0},y_{n_0}) \le \alpha$. This
  decreasing sequence will converge to a cluster point of $p(x_{n},y_{n})$
  strictly smaller than $\alpha$ which is contradictory with the definition of
  $\alpha$. Thus we have $\alpha < p(x_n,y_n)$ for all $n \in \N$ and $\alpha >0$.
  We then consider $\delta(\alpha)$ given by 
  \ref{it:acfdeux} for
  $\epsilon=\alpha$.  By definition of $\alpha$ we can find $(x_i,y_i)$ such that
  $\alpha < p(x_i, y_i) < \alpha + \delta(\alpha)$ and by 
  \ref{it:acfdeux} we will obtain $\nu \in \N$
  such that $p(x_{\nu+i},y_{\nu+i}) \le \alpha$ which contradicts $\alpha < p(x_n,y_n)$
  for all $n \in \N$. Thus we conclude that $\alpha=0$.

  We prove now that $\liminf_{n \to \infty} p(x_n,y_n)=0$ and 
  \ref{it:acfun} imply that $\limsup_{n \to \infty} p(x_n,y_n)=0$.
  For $\epsilon > 0$ given, we consider $\delta$ given by 
  \ref{it:acfun}. Since $\liminf_{n \to \infty} p(x_n,y_n)=0$ then we can
  find $i\in \N$ such that $ p(x_i, y_i) < \delta$. Thus by \ref{it:acfun} we have 
  $\limsup_{n \to \infty} p(x_{n+i},y_{n+i}) \le \epsilon$ and thus 
  successively $\limsup_{n \to \infty} p(x_{n},y_{n}) \le
  \epsilon$ and $\limsup_{n \to \infty} p(x_{n},y_{n}) =0$ and the result  follows.
\end{myproof}

\begin{defn} \label{def:acfd} We say that a sequences $\seq{x}$, with
  $x_n \in \X$ is $p$-ASF-$2$ if we have the following property:
  \begin{enumerate}[start=4]
  \item \label{it:acfddeux} For each $\epsilon > 0$, there exist $\delta > 0$ and 
    $\nu \in \N$ such that 
    if for $i$,$j \in \N$ we have $\epsilon < p(x_i, x_j) < \epsilon + \delta$, then  
    $p(x_{\nu+i},x_{\nu+j}) \le \epsilon$\,.
  \end{enumerate}
\end{defn}

Let $q$ be a given function from $\X\times \X$ into $[0,\infty)$ and $p=G\circ q$ 
where the mapping $G$ is a nondecreasing right continuous function such that 
$G(t)>0$ for $t>0$. We first show here that when a sequence is $(G\circ q)$-ASF-$2$ then 
it is also a $q$-ASF-$2$ sequence if \ref{it:acftrois-g} is satisfied by $p$. 
Note that Property~\ref{it:acfddeux} (resp. \ref{it:acftrois-g}) is a kind of uniform extension 
of \ref{it:acfdeux} (resp. \ref{it:acftrois}) when only one sequence is involved. 

\begin{lem}(\cite[In Theorem 6]{suzuki-2007})\label{lem:ASF-D} Let $\seq{x}$ be a $p$-ASF-2 
sequence and suppose that $p=G\circ q$ where $G$ is a nondecreasing right continuous 
function such that $G(t)>0$ for $t>0$. Suppose that we have 
\begin{enumerate}[start=5]
\item \label{it:acftrois-g} for each given $(x_i,x_j)$ such that $p(x_i,x_j)\ne 0$ 
  there exists $\nu \in \N$ such that $$p(x_{\nu+i}, x_{\nu+j}) < p(x_i,x_j)\,,$$
\end{enumerate}
then $\seq{x}$ is a $q$-ASF-2 sequence.
\end{lem}

\begin{myproof} The proof is contained in \cite[Theorem 6]{suzuki-2007}.
  Fix $\eta >0$ and consider $\epsilon = G(\eta)$. 
  Since $G(t) > 0$ for $t >0$ we have $\epsilon >0$. Then we can 
  use \ref{it:acfddeux} to obtain $\delta >0$ and $\nu \in \N$ such that 
  $\epsilon < p(x_i, x_j) < \epsilon + \delta$, for some $i$, $j\in \N$ 
  implies  $p(x_{\nu+i},x_{\nu+j}) \le \epsilon$\,.
  Since $G$ is nondecreasing right continuous we can find $\beta$ such that 
  $G([\eta,\eta+\beta]) \subset [\epsilon,\epsilon+\delta)$. 
  Thus suppose that $\eta < q(x_i,x_j) < \eta+\beta$, we then have 
  $\epsilon \le G(q(x_i,x_j)) < \epsilon+\delta$. Since $G$ is 
  nondecreasing it can be constant and equal to $\epsilon$ on a non 
  empty interval $[\eta,\eta+\overline{\beta}) \subset \eta+\beta$ 
  in the contrary we will have $\epsilon < G(\eta+\gamma)$ for 
  $\gamma \in(0,\beta)$. If we are in the second case then 
  $\epsilon < G(q(x_i,x_j)) < \epsilon+\delta$ and using 
\ref{it:acfddeux} we obtain 
$G(q(x_{i+\nu},x_{j+\nu})) \le \epsilon < G(\eta+\gamma)$ we 
thus have $q(x_{i+\nu},x_{j+\nu}) < \eta+\gamma$ for all 
$\gamma \in(0,\beta)$ and consequently $q(x_{i+\nu},x_{j+\nu})\le \eta$. 
In the first case we have $G(q(x_i,x_j))=\epsilon$ for
$\eta < q(x_i,x_j) < \eta+\overline{\beta}$. Using 
\ref{it:acftrois-g} we can find $\nu\in \N$ such that 
\begin{equation}
  G(q(x_{i+\nu},x_{j+\nu})) < G(q(x_i,x_j))=\epsilon = G(\eta)
\end{equation}
and thus $q(x_{i+\nu},x_{j+\nu}) \le \eta$. We thus have proved that 
Property~\ref{it:acfddeux} is satisfied by $q$.
\end{myproof}

We prove now that $p$-ASF-$2$ sequences mixed with convergence properties 
of the sequence $p(x_n,x_{n+1})$ gives $p$-Cauchy properties. More precisely we have the following lemma.

\def\qdist{r}

\begin{lem}\label{lem:ASFD} Let $\seq{x}$, be a $p$-ASF-2 sequence  
  and suppose that $p$ is such that $p(x,y) \le p(x,z)+\qdist(z,y)$ and 
  $p(x,y) \le \qdist(x,z)+ p(z,y)$ for all 
  $x$, $y$, $z\in \X$ where the mapping $\qdist:\X\times \X \to [0,\infty)$ 
  satisfies the triangle inequality $\qdist(x,y)  \le \qdist(x,z)+\qdist(z,y)$ for all 
  $x$, $y$, $z\in \X$. If the sequence $\seq{x}$ is such that 
  $\lim_{n \to \infty} \qdist(x_{n},x_{n+1})=0$ and $\lim_{n \to \infty} p(x_{n},x_{n+1})=0$ 
  then we have $\lim_{n\to \infty} \sup_{m>n} p(x_n,x_m) = 0$.
\end{lem}

\begin{myproof} 
  We follow \cite[Lemma 2]{suzuki2007} where a similar proof is given when $\qdist=p$.
  Let $\epsilon >0 $ be fixed 
  and consider $\delta$ and $\nu$ given by \ref{it:acfddeux}.
  There exists $N\in \N$ such that 
  $\qdist(x_n,x_{n+1}) < \delta/\nu$ and $p(x_n,x_{n+1}) < \delta/\nu$ for all $n \ge N$. We first have for $k\le \nu$:
  \begin{align}
    \qdist(x_{n},x_{n+k}) & \le \sum_{i=0}^{k-1}  \qdist(x_{n+i},x_{n+i+1}) < k \frac{\delta}{\nu} \le \delta
    \label{recfornu} \\
\intertext{and}
    p(x_{n},x_{n+k}) & \le p(x_{n},x_{n+1}) + \sum_{i=1}^{k-1}  \qdist(x_{n+i},x_{n+i+1}) < k \frac{\delta}{\nu} \le \delta \label{recfornd}
  \end{align}
  We suppose that $p(x_{n},x_{n+\alpha}) < \epsilon+\delta$ is satisfied for $\alpha \in[1,k]$ and 
  we want to prove that the same inequalities are satisfied for $\alpha \in[1,k+1]$. 
  Using \eqref{recfornd} we may assume that $k\ge \nu$. Using the mixed 
  triangle inequality satisfied by $p$ we have the two separate 
  inequalities:
  \begin{align}
    p(x_{n},x_{n+k+1}) & \le  p(x_{n},x_{n+k+1-\nu}) + \sum_{i=1-\nu}^{0}  \qdist(x_{n+k+i},x_{n+k+i+1})\nonumber \\
    & <  p(x_{n},x_{n+k+1-\nu})  + \delta   \label{inequn} \\
\intertext{and}   
p(x_{n},x_{n+k+1}) & \le  \qdist(x_{n},x_{n+\nu}) + p(x_{n+\nu},x_{n+k+1-\nu+\nu})
\label{ineqdeux} 
\end{align}
  By hypothesis we have $p(x_n,x_{n+k+1-\nu}) < \epsilon+\delta$. If $p(x_n,x_{n+k+1-\nu}) 
  \le \epsilon$ then using \eqref{inequn} we obtain $p(x_{n},x_{n+k+1}) < \epsilon+\delta$ 
  else we can use  \ref{it:acfddeux} to first get 
  $p(x_{n+\nu},x_{n+k+1-\nu+\nu}) \le \epsilon$ and using \eqref{recfornu} and \eqref{ineqdeux}
  we obtain $p(x_{n},x_{n+k+1}) < \epsilon+\delta$.
\end{myproof}

In \cite{suzuki-2001}, T. Suzuki introduces the definition 
of a $\tau$-distance. We just recall here two properties which are satisfied 
by $\tau$-distance: if a function $p$ 
from $\X\times \X$ into $\R^+$ is a $\tau$-distance it satisfies $p(x,y)\le p(x,z)+p(z,y)$ for all 
$x$, $y$, $z \in \X$ and if a sequence $\seq{x}$ in $\X$ satisfies 
$\lim_{n \to \infty} \sup_{m >n} p(x_{n},x_{m})=0$ then $\seq{x}$ is 
a Cauchy sequence. We thus have the following theorem.

\begin{thm}\label{thmpasf} Let $\seq{x}$ be a $p$-ASF-2 sequence in $\X$ such that 
$\seq{x}$ and $\seq{y}$ are $p$-ASF-1 for $y_{n}=x_{n+1}$ for 
all $n\in \N$. If one of the following assumptions holds true
\begin{enumerate}[labelindent=\parindent,label=(\roman*),ref=\emph{(\roman*)}]
\item \label{thmpasf:un} $p=q$ and $q$ is a $\tau$-distance\,;
\item \label{thmpasf:deux} $p=G(q)$ where $q$ is a $\tau$-distance and where $G$ is a nondecreasing right continuous function such that $G(t)>0$ for $t >0$ and \emph{\ref{it:acftrois-g}} is satisfied by the sequence $\seq{x}$ (for the mapping $p=G(q)$)\,;
\item \label{thmpasf:trois} $p$ is a $\tau$-distance such that $p(x,y) \le p(x,z)+q(z,y)$ 
  and $p(x,y) \le q(x,z)+p(z,y)$ for all 
  $x$, $y$, $z\in \X$  where the mapping $q:\X\times \X \to [0,\infty)$ 
  satisfies the triangle inequality $q(x,y) \le q(x,z)+q(z,y)$ for all 
  $x$, $y$, $z\in \X$ and $\lim_{n\to \infty} q(x_n,x_{n+1})=0$\,;
\end{enumerate}
then, $\seq{x}$ is a Cauchy sequence.
\end{thm}

\begin{myproof} First note that, in the three cases, using Lemma~\ref{lem:ASF} 
we have that $$\lim_{n\to \infty} p(x_n,x_{n+1})=0\,.$$

\ref{thmpasf:un} We consider the case $p=q$. Since $\lim_{n\to \infty} p(x_n,x_{n+1})=0$ 
We can use Lemma~\ref{lem:ASFD} (with $r=p$) to obtain $\lim_{n\to \infty} \sup_{m>n} p(x_n, x_m) = 0$ and 
since $p$ is a $\tau$-distance we obtain the fact that $\seq{x}$ is a Cauchy 
sequence \cite[Lemma 1]{suzuki-2001}. 

\ref{thmpasf:deux} Suppose now that $p=G(q)$, we have $\lim_{n\to \infty} G(q(x_n,x_{n+1}))=0$. 
This is only possible if $G(0)=0$ and we thus also obtain that 
$\lim_{n\to \infty} q(x_n,x_{n+1})=0$. Using Lemma~\ref{lem:ASF-D} we obtain that 
$\seq{x}$ is $q$-ASF-2 and we conclude as in the part \ref{thmpasf:un} using now the 
$\tau$-distance $q$. 

\ref{thmpasf:trois} Here we can use Lemma~\ref{lem:ASFD} to obtain $\lim_{n\to \infty} \sup_{m>n} p(x_n, x_m) = 0$ and using the fact that $p$ is a $\tau$-distance the conclusion follows the lines of the case \ref{thmpasf:un}.
\end{myproof}

\begin{rem} Note that we have proved during the proof of Theorem~\ref{thmpasf} 
that if we have two sequences $\seq{x}$ and $\seq{y}$ which are $p$-ASF-1 
with $p=G\circ q$ then $G(0)=0$.
\end{rem}

\section{Links with ACF sequences}
\label{defacf}
We first recall here the definition of an ACF mapping. Then we 
give a definition of a $T$-ASF mapping by defining properties 
which are to be satisfied by the sequences $\{T^nx\}$ for $x \in \X$. 
We prove in Theorem~\ref{asf-acf} that the two definitions are equivalent. 

\begin{defn}\cite[Definition 1]{suzuki2007} Let $(\X, d)$ be a metric space. 
  Then a mapping $T$ on $\X$ is said to be an asymptotic
  contraction of the final type (ACF, for short) if the following hold:
\begin{enumerate}[label=\emph{($\mbox{D}_{\arabic*}$)}]
\item $\lim_{\delta \to 0^+} \sup \left\{ \limsup_{n\to \infty} d(T^n x,T^n y): d(x, y) <\delta
\right\} = 0$.
\item  For each $\epsilon > 0$, there exists $\delta > 0$ such that for $x$, $y \in \X$
  with $\epsilon < d(x, y) < \epsilon + \delta$, there exists $\nu \in N$ such that 
  $d(T^\nu x,T^\nu y) \le \epsilon$.
\item  For $x$, $y \in \X$ with $x\not=y$, there exists $\nu \in N$ such that 
  $d(T^\nu x,T^\nu y)<d(x,y)$.
\item  For $x \in \X$ and $\epsilon >0$, there exist $\delta >0$ and $\nu \in N$ such that
\begin{equation}
  \epsilon <d(T^i x,T^j x) < \epsilon +\delta \;\mbox{implies}\; 
  d(T^\nu \circ T^i x,T^\nu\circ T^j x )\le \epsilon\,
\end{equation} 
for all $i$, $j \in N$.
\end{enumerate}
\end{defn}

\begin{thm}\label{asf-acf} Let $(\X, d)$ be a metric space. A mapping $T$ on $\X$ is 
  said to be $p$-ASF if for all $x$, $y \in \X$ the sequences 
  $\{T^n x\}$ and   $\{T^n y\}$ are $p$-ASF-1 and $\{T^n x\}$ is $p$-ASF-2. 
  Then, $T$ is an $ACF$ mapping is equivalent to $T$ is a $d$-ASF mapping. 
\end{thm}

\begin{myproof}
Suppose that the mapping $T$ is $ACF$. For each $x$, $y\in \X$ it is very easy to check 
and left to the reader that $\{T^n x\}$ and $\{T^n y\}$ are $d$-ASF-1 and 
$\{T^n x\}$ is $d$-ASF-2. Thus, $T$ is $d$-ASF. 

If $T$ is $d$-ASF, using Lemma~\ref{lem:ASF} we obtain 
$\lim_{n\to \infty} d(T^n x,T^n y) = 0$. 
If we consider the special case $y=Tx$ and the sequence $x_n = T^nx$ we 
obtain using Theorem~\ref{thmpasf} that $\{T^nx\}$ is a Cauchy sequence. Then 
using \cite[Theorem 6]{suzuki2007}\footnote{
 We first recall from \cite[Theorem 6]{suzuki2007} that for a mapping 
$T$ on a metric space $(\X, d)$ the following are equivalent: 
\begin{enumerate}[labelindent=\parindent,label=(\roman*),ref=\emph{(\roman*)}]
\item $T$ is an ACF.
\item $\lim_{n\to \infty} d(T^n x,T^n y) = 0$ holds true and $\{T^nx\}$ is a Cauchy sequence 
  for all $x$, $y \in \X$.
\end{enumerate}
}
we obtain that the mapping $T$ is ACF. 
\end{myproof}

Existence and uniqueness of fixed points of $p$-ASF mappings is now
obtained. Note that, in the special case where the mapping $p$ is equal to $d$
(i.e when we use the $\tau$-distance $p=d$ in $(i)$) the next theorem gives same
results as \cite[Theorem 5]{suzuki2007}.

\begin{thm}\label{fixedpt} Let $(\X,d)$ be a complete metric space, $T$ be a $p$-ASF mapping which is such that 
$T^l$ is continuous for some $l\in \N$ ($l>0$). 
We suppose that the function $q$ is a $\tau$-distance and 
one of the following holds true for the mapping $p$:
\begin{enumerate}[labelindent=\parindent,label=(\roman*),ref=\emph{(\roman*)}]
\item $p=q$\,.
\item $p=G(q)$ where $G$ is a nondecreasing right continuous function such that $G(t)>0$ for $t >0$ and \emph{\ref{it:acftrois-g}} is satisfied by the sequence $\seq{x}$ (for the mapping $p=G(q)$). 
\end{enumerate}
then, there exists a fixed point $z\in\X$ of $T$. Moreover, if for every sequences $\seq{x}$ and 
$\seq{y}$ $\lim_{n\to \infty} p(x_n,y_n)=0$ implies that $\lim_{n\to \infty} d(x_n,y_n)=0$ then the 
fixed point is unique and $\lim_{n\to \infty} T^nx=z$ holds true for every $x\in \X$
\end{thm}

\begin{myproof} For every $x\in \X$ using Theorem~\ref{thmpasf} we know that $\{T^nx\}$ is a 
Cauchy sequence. By Lemma~\ref{lem:ASF} we know that $\lim_{n\to \infty} p(T^nx,T^ny)=0$. 
We then have all the ingredients of  \cite[Theorem 4 and Lemma 3]{suzuki2007} to conclude the proof.
\end{myproof}

\section{An application to $ACF$ cyclic mappings} 
\label{cycling}

We suppose here that $\X$ is a uniformly convex Banach space and thus 
$d(x,y)\defpar \norm{x-y}$. We consider $A$ and $B$ two nonempty subsets of 
$\X$, $A$ being convex, and a cyclic mapping $T: A \cup B \to A \cup B$. 
We recall that $T$ is a cyclic mapping if $T(A) \subset B$ and $T(B) \subset A$. 
We define a mapping $p: \X\times \X \to \R^+$ 
by $p(x,y) \defpar d(x,y) -d(A,B)$ where $d(A,B) \defpar \inf \{ d(x,y) \,|\, x\in A, y\in B\}$. 

Then, using previous results we can give a short proof of a theorem 
which extends~\cite[Theorem 1]{suzuki-2008}.

\begin{thm}\label{thm:bcauchy} Suppose that the mapping $T$ is $p$-ASF, then the sequence $\{T^{2n}x\}$ 
for $x\in A$ is a $d$-Cauchy sequence. 
\end{thm}

\begin{myproof} For a given $x\in \X$, we consider the sequence $x_n=T^nx$.
  Since $T$ is $p$-ASF we have by Lemma~\ref{lem:ASF} that $\lim_{n\to \infty}
  p(x_n,x_{n+1})=0$.  Using the definition of $p$ we immediately also have
  $\lim_{n\to \infty} d(x_n,x_{n+1})=d(A,B)$. Using Lemma~\cite[Lemma
  4]{suzuki-2008} we obtain that $\lim_{n \to \infty} d(x_{2n},x_{2n+2})=0$
  (convexity of $A$ and uniformly convexity of $\X$ are used here). We now
  consider the sequence $\{T^{2n} x\}$ taking values in $A$.  We have
  $\lim_{n\to \infty} p(x_{2n},x_{2n+2})=0$ and as it was already shown $\lim_{n
    \to \infty} d(x_{2n},x_{2n+2})=0$. If the sequence $x_n=T^nx$ is $p$-ASF-2
  then it is the same for the sequence $\{T^{2n} x\}$. The distance $d$ satisfy
  the triangle inequality and it is straightforward to see that we have the two
  mixed triangle inequality $p(x,y) \le p(x,z)+d(z,y)$ and $p(x,y) \le
  d(x,y)+p(z,y)$ for all $x$, $y$, $z\in \X$. We can thus apply
  Lemma~\ref{lem:ASFD} to the sequence $\{T^{2n} x\}$ with $x\in A$ to obtain
  that it is a $p$-Cauchy sequence. It is now easy to see by contradiction that
  a $p$-Cauchy sequence is a $d$-Cauchy sequence \cite[Proof of Theorem
  2]{suzuki-2008}.  The key argument being again the use of \cite[Lemma
  4]{suzuki-2008}
\end{myproof}

We extend now \cite[Theorem 2]{suzuki-2008} which was stated for continuous
cyclic Meir-Keeler contractions to continuous $p$-ASF mappings.

\begin{thm}\label{thm:cyclic}Suppose in addition that $A$ is closed, 
  $T$ is $p$-ASF and $T^l$ is continuous for
  some $l\in\N$ ($l>0$) then there exists a unique best 
  proximity point $z\in A$ (i.e
  $d(z,Tz)=d(A,B)$). Moreover $\lim_{n\to \infty} T^{2n} x=z$ for each $x\in A$.
\end{thm}

\begin{myproof} Using Theorem~\ref{thm:bcauchy}, the sequence $\{T^{2n}x\}$ for each $x\in A$ is a $d$-Cauchy sequence. Using Lemma~\ref{lem:ASF}, we have $\lim_{n \to \infty} p(T^{n}x,T^{n}y)=0$ for each $x$, $y \in A$, hence for $(x,Tx)$ it gives $\lim_{n \to \infty} p(T^{2n}x,T^{2n+1}x)=0$ and for $(Tx,y)$ it gives $\lim_{n \to \infty} p(T^{2n+1}x,T^{2n}y)=0$. Using again \cite[Lemma 4]{suzuki-2008} we obtain 
$\lim d(T^{2n}x,T^{2n}y)=0$ and we can use  \cite[Theorem 4 and Lemma 3]{suzuki2007} to conclude the proof.
\end{myproof}

\section{ASMK Sequences}

We introduce in this section the definition of ASMK sequences. 
It is an adaptation to sequences of the ACMK (Asymptotic contraction of 
Meir-Keeler type) definition used for mappings \cite{suzuki-2006}. 
It is proved in \cite[Theorem 3]{suzuki2007} that an ACMK mapping 
on a metric space is an ACF mapping. We will prove in this section 
similar results which relate ASMK sequences to ASF sequences. These 
results will be used in next section for studying sequences of 
alternating mappings. 

\begin{defn} \label{Famk} We say that two sequences $\seq{x}$, $\seq{y}$ with
  $x_n,y_n \in \X$ are $p$-ASMK-$1$ if there exists a sequence $\{\psi_n\}$ 
  of functions from $[0,\infty)$ into itself satisfying $\psi_n(0)=0$\footnote{Note that this assumption can be removed when $F(0)>0$.}
  for all $n \in \N$ and the following:
  \begin{enumerate}[start=6]
  \item \label{it:famkun} $\limsup_n \psi_n(\epsilon) <  \epsilon \texte{for all} \epsilon > 0$.
  \item \label{it:famkdeux}  For each $\epsilon >0$, there exists $\delta >0$
    such that for each $t \in [\epsilon, \epsilon + \delta]$ there exists 
    $\nu \in \N$ such that $\psi_\nu (t) < \epsilon$.
  \item \label{it:famktrois} $F\bp{p(x_{n+i},y_{n+i})} \le \psi_n\Bp{F\bp{p(x_i,y_i)}}$ for
    all $n$, $i \in \N$. $F$ is a given right continuous nondecreasing mapping
    such that $F(t) >0$ for $t\not=0$.
  \end{enumerate}
\end{defn}

\begin{lem}\label{lem:asmk-asf} Suppose that the two sequences $\seq{x}$, $\seq{y}$ are $p$-ASMK-$1$ then they are $p$-ASF-1.  
\end{lem}

\begin{myproof} 
  \ref{it:acfun}: For all $n$, $i \in \N$ we have by~\ref{it:famktrois} 
  and~\ref{it:famkun} when $F\bp{p(x_i,y_i}\ne 0$ that 
  \begin{align}
    F\bp{p(x_{n+i},y_{n+i})} & \le \psi_n\Bp{F\bp{p(x_i,y_i)}} \nonumber \\
    & \le \limsup_{n\to \infty} \psi_n\Bp{F\bp{p(x_i,y_i)}}  \nonumber \\
    & < F\bp{p(x_i,y_i)}\,. \nonumber 
  \end{align}
  Since $F$ is nondecreasing and the inequality is strict we obtain for all 
  $n \in \N$:
  \begin{align*}
    p(x_{n+i},y_{n+i}) < p(x_i,y_i) \,,
  \end{align*}
  and thus 
  \begin{align*}
    \limsup_{n\to \infty}  p(x_{n+i},y_{n+i})  \le p(x_i,y_i) \,.
  \end{align*}
  Then \ref{it:acfun} follows easily when $F\bp{p(x_i,y_i}\ne 0$. When $F\bp{p(x_i,y_i}= 0$, 
  we have by~\ref{it:famktrois}  $F\bp{p(x_{n+i},y_{n+i})}\le 0$ for all $n\in \N$. 
  Since $F$ is a right continuous mapping such that $F(t) >0$ we must have $F(0)\ge 0$. 
  Thus we have that $F\bp{p(x_{n+i},y_{n+i})}= 0$ for all $n\in \N$ and thus
  $p(x_{n+i},y_{n+i})=0$ for all $n\in \N$ and the same conclusion holds. 
  \ref{it:acfdeux}:
  for $\epsilon > 0$ we know that $F(\epsilon) >0$ and we can use
  \ref{it:famkdeux} to find $\delta >0$ such that for each 
  $t \in [F(\epsilon), F(\epsilon) + \delta]$ we can find $\nu \in \N$ such that $\psi_\nu 
  (t) < F(\epsilon)$. Since
  $F$ is right continuous and nondecreasing we can find $\delta'$ such that
  $F([\epsilon,\epsilon+\delta']) \subset [F(\epsilon), F(\epsilon) + \delta]$.
  Thus, taking $i\in \N$ such that $\epsilon < p(x_i,y_i) < \epsilon + \delta'$, 
  we can find $\nu$ such that $\psi_{\nu}(F(p(x_i,y_i))) < F(\epsilon)$. And we 
  conclude using \ref{it:famktrois} that:
  \begin{equation}
    F\bp{p(x_{\nu+i},y_{\nu+i})} \le \psi_{\nu} \Bp{F\bp{p(x_i,y_i)}} < F(\epsilon) \le F\bp{p(x_i,y_i)}\,.
  \end{equation}
  Thus we have $F\bp{p(x_{\nu+i},y_{\nu+i})}  <  F\bp{p(x_i,y_i)}$ and since $F$ is nondecreasing 
  and the inequality is strict we obtain \ref{it:acfdeux}. 
  
  \ref{it:acftrois}: Let $i$ be given such that $p(x_i,y_i)\not=0$ and start as in the previous 
  paragraph using $\epsilon =p(x_i,y_i)$. We can find $\nu \in \N$ such that 
  $\psi_{\nu}\Bp{F\bp{p(x_i,y_i)}} < F(\epsilon)$ which combined with \ref{it:famktrois} gives:
  \begin{equation}
    F\bp{p(x_{\nu+i},y_{\nu+i})} \le \psi_{\nu} \Bp{F\bp{p(x_i,y_i)}} < F(\epsilon) = F\bp{p(x_i,y_i)}\,.
  \end{equation}
  Since $F$ is nondecreasing and the inequality is strict the result follows.
\end{myproof}

\begin{defn} We say that two sequences $\seq{x}$, $\seq{y}$ with
  $x_n,y_n \in \X$ are $p$-ASMK-$2$ when \emph{\ref{it:famktrois}} is replaced by 
  \begin{enumerate}[start=9]
  \item \label{it:famktrois-p} $F\bp{p(x_{n+i},y_{n+j})} \le \psi_n \Bp{F\bp{p(x_i,y_j)}}$ 
    for all $n$, $p$, $i$,$j \in \N$.
  \end{enumerate}
\end{defn}

\begin{cor} If two sequences $\seq{x}$, $\seq{y}$ with $y_{n}=x_{n+1}$ are $p$-ASMK-$2$ then they 
are $p$-ASF-1 and the sequence $\seq{x}$ is $p$-ASF-2. Moreover, assumption \emph{\ref{it:acftrois-g}} holds true 
for $p$.
\end{cor}

\begin{myproof} It is obvious to see that if two sequences $\seq{x}$, $\seq{y}$ are 
   $p$-ASMK-$2$ then they are $p$-ASMK-$1$. Thus by Lemma~\ref{lem:asmk-asf} they are 
   $p$-ASF-1. Proving that \ref{it:acfddeux} holds true is similar to the proof that \ref{it:acfdeux} 
   holds true in Lemma~\ref{lem:asmk-asf} and proving that \ref{it:acftrois-g} holds true 
   follows the same steps as the proof that \ref{it:acftrois} holds true in Lemma~\ref{lem:asmk-asf}.
\end{myproof}

\section{A sequence of alternating mappings}
\label{alternate}
In this section $p$ is a given function from $\X\times \X$ into $[0,\infty)$ 
such that such $p(x,y) \le p(x,z)+p(z,y)$ for all $x$, $y$, $z\in \X$ and 
$p(x,y)=p(y,x)$ for all $x$, $y \in \X$.

\begin{defn} We will say that the pair $(T,S)$ satisfy the $(F,\psi)$-contraction 
property if we can find two functions $F$ and $\psi$ such that:
\begin{equation}
  F \bp{ p(Tx,Sy)} \le \psi \Bp{F\bp{ M(x,y)}}
\end{equation}
where 
\begin{equation}
  M(x,y) \defpar \max \left\{p(x,y),p(Tx,x),p(Sy,y), 
    \frac{1}{2}\left\{p(Tx,y)+p(Sy,x) \right\} \right\}
\end{equation}
The function $F:\R^+ \to \R^+$ is a given right continuous nondecreasing mapping
such that $F(t) >0$ for $t\not=0$. The function $\psi:\R^+ \to \R^+$ is a given 
nondecreasing upper semicontinuous function satisfying $\psi(t)< t$ for each $t>0$ and 
$\psi(0)= 0$. 
\end{defn} 

We first start by a technical lemma.

\begin{lem}\label{lemme_dix} Let the pair of mappings $(T,S)$ be a 
$(F,\psi)$-contraction. 
Suppose that $x=S\alpha$ and $p(x,Tx)\not=0$ then we have:
\begin{equation}
  \label{ineqT}
  F \bp{ p(x,Tx)} \le \psi\Bp{ F \bp{p(S\alpha,\alpha)}} \,.
\end{equation}
Suppose that $y=T\alpha$ and $p(y,Sy)\not=0$ then we have:
\begin{equation}
  \label{ineqS}
  F \bp{p(Sy,y)} \le \psi \Bp{ F\bp{p(\alpha,T\alpha)}} \,.
\end{equation}
\end{lem}

\begin{myproof}
We prove the first inequality \eqref{ineqT}. We suppose that $x=S\alpha$ then 
we have
\begin{align}
  F \bp{p(x,Tx)} &= F \bp{p(Tx,x)} =  F \bp{p(Tx,S\alpha)} \le \psi \Bp{ F \bp{M(x,\alpha)}} \nonumber 
\end{align}
and we have:
\begin{align}
  M(x,\alpha) &=  \max \left\{p(x,\alpha),p(Tx,x),p(S\alpha,\alpha), 
    \frac{1}{2}\left\{p(Tx,\alpha)+p(S\alpha,x) \right\} \right\}  \nonumber \\
  &= \max \left\{p(x,\alpha),p(Tx,x),\frac{1}{2} p(Tx,\alpha) \right\}  \nonumber \\
  &= \max \left\{ p(x,\alpha),p(Tx,x) \right\}\,.
\end{align}
We show now that the maximum cannot be achieved by $p(Tx,x)$. Indeed, 
suppose that $M(x,\alpha)=p(Tx,x)$ then we would have 
\begin{align}
  F \bp{p(x,Tx)} & \le  \psi \Bp{ F \bp{p(x,Tx)}} \nonumber 
\end{align}
which is not possible since $p(x,Tx)\ne 0$ and there does not exist 
$x > 0$ such that $F(x) \le \psi(F(x))$ (since for $x>0$ we have that 
$F(x) \le \psi(F(x)) < F(x)$).
The proof for the second inequality is very similar and thus omited. 
\end{myproof}

We now introduce the alternating 
sequence of mappings $\seq{\Gamma}$ defined by:
\begin{equation}
  \Gamma_n \defpar \begin{cases} T, & \mbox{if } n\mbox{ is even} \\ S, & \mbox{if } n\mbox{ is odd}\,. \end{cases}
\end{equation}
Then, we consider the two sequences $\seq{x}$ and $\seq{y}$ defined by 
\begin{equation}
x_{n+1} =\Gamma_{n} x_n \quad\mbox{and} \quad y_{n+1}= \Gamma_{n+1} y_n\,.
\label{defxy}
\end{equation} 
It is very easy to check that when the two sequences 
are initiated with $(x_0,y_0)=(Sx,x)$ for a given $x\in \X$ they are 
related by $y_{n+1}=x_n$ and that only the two following cases can occur:
\begin{equation}
  (x_{n+1},y_{n+1})= \begin{cases} (Sx_n,x_n)& \text{with }\, x_n=Tx_{n-1} 
  \text{ and }\,  y_n = x_{n-1}   \\
  (T x_n,x_n) & \text{with }\, x_n= S x_{n-1} \text{ and }\, y_n = x_{n-1} \end{cases}
\end{equation}
If we are in the first (resp. the second) case we use \eqref{ineqS} (resp. 
\eqref{ineqT}) to obtain the inequality
\begin{equation}
  F (p(x_{n+1},y_{n+1})) \le \psi (F(p(x_n,y_n)))\label{ineqfpu}
\end{equation}
We thus have the following easy lemma:

\begin{lem}\label{lem:tscontr}Let the pair of mappings $(T,S)$ be a 
$(F,\psi)$-contraction and $\seq{x}$ and $\seq{y}$ be two 
sequences defined by \eqref{defxy}. If the two sequences 
are initiated by $(x_0,y_0)=(Sx,x)$ we have $y_{n+1}=x_n$ and
\begin{equation}
  F \bp{p(x_{n+1},y_{n+1})} \le \psi \Bp{ F \bp{ p(x_n,y_n)}} \label{ineqfp}\,.
\end{equation}
If $x_{n+1}=Tx_n$ (resp. $x_{n+1}=Sx_n$) and $y_{k+1}= Sy_k$ (resp. $y_{k+1}= Ty_k$)
we have:
\begin{equation}
  F \bp{p(x_{n+1},y_{k+1})} \le \psi \Bp{ F \bp{p(x_n,y_k)+\max(p(x_{n},x_{n+1}),p(x_{k},x_{k+1}))}}\,.
  \label{ineqmelange}
\end{equation}
\end{lem}

\begin{myproof} Since inequation \eqref{ineqfp} was proved by \eqref{ineqfpu},
it just remains to prove inequality \eqref{ineqmelange}. 
Suppose that $x_{n+1}=Tx_n$ and $y_{k+1}= Sy_k$ then we have
\begin{align}
  M(x_n,y_k) &=\max \left\{p(x_n,y_k),p(Tx_n,x_n),p(Sy_k,y_k), 
    \frac{1}{2}\left\{p(Tx_n,y_k)+p(Sy_k,x_n) \right\} \right\} \nonumber \\
  &= \max \left\{p(x_n,y_k),p(x_{n+1},x_n),p(y_{k+1},y_k), 
    \frac{1}{2}\left\{p(x_{n+1},y_k)+p(y_{k+1},x_n) \right\} \right\} \nonumber \\
  &\le \max \left\{p(x_n,y_k),p(x_{n+1},x_n),p(y_{k+1},y_k),\right. \nonumber \\
  &\hspace{1.5cm} \left. p(x_{n},y_k)  + \max(p(x_{n+1},x_n),p(y_{k+1},y_k)) \right\}  \nonumber \\
  &\le  p(x_{n},y_k)+ \max(p(x_{n+1},x_n),p(y_{k+1},y_k)) 
\end{align}
We thus have 
\begin{align}
  F \bp{p(x_{n+1},y_{k+1})} &\le \psi \Bp{F \bp{M(x_n,y_k)}}  \nonumber \\
  & \le \psi \Bp{F \bp{ p(x_{n},y_k)+\max(p(x_{n+1},x_n),p(y_{k+1},y_k))}} \,. \label{ineqcauchy}
\end{align}
If the opposite situation is $x_{n+1}=Sx_n$ and $y_{k+1}=Ty_k$ we obtain 
the same result by the same arguments. 
\end{myproof}

We make here a direct proof of the fact that the 
sequence $\seq{x}$ is a Cauchy sequence when $\lim_{n\to \infty} p(x_{n+1},x_n)=0$ 
is assumed. This last property will be derived from $p$-ASMK-$1$ properties as proved in 
Theorem~\ref{thm:ptfixe}.

\begin{lem} Let the pair of mappings $(T,S)$ be a $(F,\psi)$-contraction. 
Suppose that $$\lim_{n\to \infty} p(x_{n+1},x_n)=0\,,$$ then the sequence $\seq{x}$ given by~\eqref{defxy} 
is a Cauchy sequence.
\label{cauchylem}
\end{lem}

\begin{myproof}
We follow here \cite{zhang2007} to prove the result by contradiction.

If the sequence is not a Cauchy sequence we can find two subsequences
$\sigma(n)$ and $\rho(n)$ such that for all $n \in \N$ 
$p(x_{\sigma(n)}, x_{\rho(n)}) \ge 2\epsilon$ and $\sigma(n) < \rho(n)$. 
Since the sequence $\{ p(x_n,x_{n+1})\}_{n\in\N}$ converges to zero we can choose $N$ 
such that $p(x_n,x_{n+1}) < \epsilon$ for all $n \ge N$. Using the triangle inequality 
$$
p(x_{\sigma(n)},x_{\rho(n)+1}) \ge p(x_{\sigma(n)},x_{\rho(n)}) - p(x_{\rho(n)+1},x_{\rho(n)})\, 
$$
we obtain that $p(x_{\sigma(n)},x_{\rho(n)+1}) > \epsilon$ for large $n$. 
Thus, we can always change the subsequence ${\rho(n)}$ in such a way that the parity between 
$\sigma(n)$ and $\rho(n)$ is conform to the one we need for 
applying inequality \eqref{ineqcauchy} and such that for all $n \in \N$
$p(x_{\sigma(n)}, x_{\rho(n)}) > \epsilon$. 

We now define $k(n)$ as follows:
\begin{equation}
  k(n) \defpar \min \left\{ 
    k > \sigma(n) \, \vert \, p(x_{\sigma(n)},x_{k}) > \epsilon \quad 
    \mbox{ with same parity as $\rho(n)$} 
  \right\}
\end{equation}
$k(n)$ is well defined and by construction $\sigma(n) < k(n)\le \rho(n)$. 
We now have that:
\begin{align}
  \epsilon < p(x_{\sigma(n)},x_{k(n)}) \le p(x_{\sigma(n)},x_{k(n)-2}) + 
  p(x_{k(n)-2},x_{k(n)}) \le \epsilon + p(x_{k(n)-2},x_{k(n)})\,.
\end{align}
the sequence $\{p(x_{k(n)-2},x_{k(n)})\}_{n\in\N}$ converges to zero since we have 
$$p(x_{k(n)-2},x_{k(n)})\le p(x_{k(n)-2},x_{k(n)-1}) + p(x_{k(n)-1},x_{k(n)}) $$
and thus $p(x_{\sigma(n)},x_{k(n)}) \to \epsilon^+$ when $n$ goes to infinity. We also 
obtain that $p(x_{\sigma(n)-1},x_{k(n)-1}) \to \epsilon$ when $n$ goes to infinity since:
\begin{align}
  |p(x_{\sigma(n)},x_{k(n)}) - p(x_{\sigma(n)-1},x_{k(n)-1})| 
  \le p(x_{k(n)},x_{k(n)-1}) + p(x_{\sigma(n)},x_{\sigma(n)-1})\,.
\end{align}
We now use inequality~\eqref{ineqcauchy} to obtain 
\begin{align}
  F \bp{p(x_{\sigma(n)},x_{k(n)})} \le \psi \Bp{F\bp{p(x_{\sigma(n)-1},x_{k(n)-1})+\delta_n}} 
\end{align}
where $\delta_n \defpar \max \bp{p(x_{\sigma(n)-1},x_{\sigma(n)}),p(x_{k(n)-1},x_{k(n)})}$. 
When $n$ goes to infinity, using the facts that $F$ is right continuous and nondecreasing 
and $\psi \circ F$ is upper semicontinuous  we obtain that $F(\epsilon) \le \psi \bp{F(\epsilon)}$ 
which is a contradiction. 
\end{myproof}

\begin{rem}\label{rem:zhang} The proof remains valid is we assume as in \cite{zhang2007} 
that the function $F$ is nondecreasing and continuous with $F(0)=0$ and $F(t)> 0$ for $t>0$ 
and that the function $\psi:\R^+ \to \R^+$ is assumed to be 
nondecreasing and right upper semicontinuous and satisfy $\psi(t)< t$ for each $t>0 $.  
The idea is to build the sequences choosing the parity so as to use~\eqref{ineqcauchy} 
in the reverse situation where
\[F \bp{p(x_{\sigma(n)+1},x_{k(n)+1})} \le \psi \Bp{F\bp{p(x_{\sigma(n)},x_{k(n)})+\delta'_n}}, \]
and where $\delta'_n \defpar \max \bp{p(x_{\sigma(n)+1},x_{\sigma(n)}),p(x_{k(n)+1},x_{k(n)})}$. 
\end{rem}

\begin{thm}\label{thm:ptfixe} Consider two mappings $T: \X \to \X $ and $S: \X \to \X$ and 
suppose that the pair $(T,S)$ has the $(F,\psi)$-contraction 
property. Let the sequence of function $\seq{\psi}$ be defined by 
$\psi_n \defpar \overset{n}{\overbrace{\psi\circ\psi\circ\cdots\psi}}$ and 
assume that \emph{\ref{it:famkun}} and \emph{\ref{it:famkdeux}} are satisfied, 
then the sequence $\seq{x}$ defined by \eqref{defxy} and initialized by $x_0=Sx$ 
is a Cauchy sequence. 
\end{thm}

\begin{myproof} The only point to prove is that assumption~\ref{it:famktrois} is satisfied. 
We consider the sequence $\seq{x}$ and the sequence $\seq{y}$ defined by \eqref{defxy} and 
initialized by $y_0=x$. Using the fact that $\psi$ is non-decreasing, we 
repeatedly use Equation~\eqref{ineqfp} in Lemma~\ref{lem:tscontr} 
to obtain assumption \ref{it:famktrois} and conclude that the two sequences 
$\seq{x}$ and $\seq{y}$ are $p$-ASMK-$1$ and then by Lemma~\ref{lem:asmk-asf} and 
\ref{lem:ASF} we obtain that $\limsup_{n \to \infty} p(x_n,x_{n+1})=0$. 
Using Lemma~\ref{cauchylem} we conclude that $\seq{x}$ is a Cauchy sequence.
\end{myproof}

We make a link here with the result of \cite{zhang2007} where
it is assumed that $F(0)=0$ and $F(t)> 0$ for $t>0$ and 
$F$ is supposed to be nondecreasing and continuous. 
The function $\psi:\R^+ \to \R^+$ is assumed to be 
nondecreasing and right upper semicontinuous and satisfy $\psi(t)< t$ for each $t>0 $ 
and $\lim_{n\to \infty} \psi_n(t)=0$. It is proved in \cite{zhang2007} that $F(x) \le \psi(F(x))$ implies 
$x=0$. We prove in the next lemma that these properties of functions $F$ and $\psi$ imply 
Properties~\ref{it:famkun} and~\ref{it:famkdeux}. 

\begin{lem} Let $\psi:\R^+ \to \R^+$ be a nondecreasing, right upper semicontinuous 
  function satisfying $\psi(t)< t$ for each $t>0$. Then the sequence of functions 
  $\seq{\psi}$ defined by 
  $\psi_n \defpar \overset{n}{\overbrace{\psi\circ\psi\circ\cdots\psi}}$ 
  satisfy \emph{\ref{it:famkun}} and \emph{\ref{it:famkdeux}}.
\end{lem}

\begin{myproof} \ref{it:famkun}: For $t>0$, since $\psi$ is nondecreasing and $\psi(t) < t$ we have $\psi_n(t) \le \psi(t) < t$ and thus \ref{it:famkun} follows. 
\ref{it:famkdeux}: Using \cite[Theorem 2]{jacek-1997} we can find a right continuous function $\overline{\psi}: \R^+ \to \R^+$ such that $\psi(t)\le \overline{\psi}(t) <t$ 
for $t>0$. Thus we easily have \ref{it:famkdeux}, since proving \ref{it:famkdeux} (using $\nu=1$) for a right continuous function is easy. 
\end{myproof}

In \cite{zhang2007} it is proved that $T$ and $S$ have a common fixed point when 
$\X$ is a complete metric space and $p=d$. The proof follows the following steps:
Since $\seq{x}$ is a Cauchy sequence it converges to $\overline{x} \in \X$. 
Using the definition of $M$ one easily checks that $M(x_{2n},\overline{x})\to d(S\overline{x},\overline{x})$ and 
$M(x_{2n},\overline{x})\ge d(S\overline{x},\overline{x})$. Moreover $Tx_{2n}=x_{2n+1}$ also converges to $\overline{x}$. We therefore have 
\begin{equation} 
  F \bp{d(Tx_{2n},S\overline{x})} \le \psi \Bp{F \bp{M(x_{2n},\overline{x})}} \,.
\end{equation}
Using next Lemma~\ref{lem:ptfixe} we obtain that $\overline{x}=S\overline{x}$. 
Then proving that $\overline{x}$ is also a fixed point of $T$ is given in 
\cite[Theorem 1]{zhang2007}. 
We therefore conclude that in order to obtain convergence of the sequence $\seq{x}$ 
to the unique fixed point of $T$ and $S$ requires to add continuity of $F$ in the hypothesis 
of Theorem~\ref{thm:ptfixe}.

\begin{lem}\label{lem:ptfixe} Suppose that $F$ is a continuous nondecreasing function, 
$\psi$ is a right upper semicontinuous function satisfying one of 
the following property:
\begin{enumerate}[label=\emph{($\mbox{E}_{\arabic*}$)}]
\item $\psi(t) < t$ for all $t >0$;
\item $\psi$ is nondecreasing and for each $t>0$, there exists 
  $\nu \in \N$, $\nu \ge 1$ such that $\psi_{\nu}(t) < t$. 
\end{enumerate}
Suppose that we have two sequences 
$\seq{\alpha}$ and $\seq{\beta}$ such that:
\begin{equation}
  F(\alpha_n) \le \psi \bp{F(\beta_n)} \,.
\end{equation}
If  $\lim_{n\to \infty} \alpha_n = \lim_{n\to \infty} \beta_n =\gamma$ and $\beta_n\ge \gamma$ for all 
$n \in \N$ then we must have $\gamma=0$.
\end{lem}
\begin{myproof} We have 
  \begin{align}
    F(\gamma)=\lim_{n\to \infty} F(\alpha_n) \le \limsup_{n\to \infty} \psi \bp{F(\beta_n)}
    \le \psi \bp{ \limsup_{n\to \infty} F(\beta_n)} \le \psi \bp{F(\gamma)};
  \end{align}
  If $\gamma \not=0$ and $\psi(\gamma) < \gamma$ we conclude that $F(\gamma) < F(\gamma)$ which is a 
  contradiction. 
  If $\psi$ is nondecreasing we consider the value of $\nu$ associated to $\gamma$ to 
  obtain:
  $F(\gamma) \le \psi_{\nu} \bp{F(\gamma)} < F(\gamma)$ and conclude again by 
  contradiction. 
\end{myproof}

\begin{rem} Note that using \cite[Theorem 2]{jacek-1997} we obtain that 
  a right upper semicontinuous function $\psi$ satisfying 
  $\psi(t)<t$ for all $t >0$ satisfies Property~\ref{it:famkdeux}
  with $\nu=1$.
\end{rem}

\bibliographystyle{elsart-num-sort.bst}
\bibliography{meir-keeler}

\end{document}